\documentclass[12pt]{article}
\usepackage{amssymb}
\def\Bbb#1{\mathbb#1}
\newtheorem{theorem}{Theorem}

\newtheorem{lemma}[theorem]{Lemma}
\newtheorem{proposition}[theorem]{Proposition}

\def\proof{\noindent{\it Proof: }}

\newcommand{\dive}{\mbox{div}}
\newcommand{\po}{{\hspace*{-1ex}}{\bf .  }}
\newcommand{\C}{\mathbb{C}}
\newcommand{\R}{\mathbb{R}}

\newcommand{\He}{{\cal H}}
\newcommand{\nab}{\bar\nabla}
\def\bea{\begin{eqnarray*} }
\def\eea{\end{eqnarray*} }
\def\be{\begin{equation} }
\def\ee{\end{equation} }
\def\<{\langle }
\def\>{\rangle }

\def\qed{\ifhmode\unskip\nobreak\fi\ifmmode\ifinner\else
\hskip5 pt \fi\fi\hbox{\hskip5 pt \vrule width4 pt
height6 pt  depth1.5 pt \hskip 1pt }}
\begin{document}

\title{\LARGE The Dirichlet problem for CMC surfaces\\ in Heisenberg space
 }\author{{ L. J. Al\'\i
as\footnote{Partially supported by MEC/FEDER project MTM2004-04934-C04-02
and Fundaci\'{o}n S\'{e}neca project 00625/PI/04, Spain.} , M. Dajczer and H.
Rosenberg}}
\date{}
\maketitle

\begin{quote}
We study constant mean curvature graphs in the Riemannian 3-dimensional Heisenberg spaces $\He=\He(\tau)$. Each such $\He$ is the total space of a Riemannian submersion onto the Euclidean plane $\R^2$ with geodesic fibers the orbits of a Killing field.
We prove the existence and uniqueness of CMC graphs in $\He$ with respect to the Riemannian submersion
over certain domains $\Omega\subset\R^2$ taking on prescribed boundary values.
\end{quote}

\section {Introduction}  In recent years, there has been much research on minimal and constant mean curvature surfaces (CMC) in the simply connected homogeneous 3-manifolds, other than space forms. Figueroa,  Mercuri and Pedrosa \cite{fmp}
gave many interesting such surfaces in $\He$, each invariant by Killing
vector fields of the ambient space. Daniel \cite{da} and Abresch-Rosenberg
\cite{ar}, \cite{ar2} have also obtained some interesting results on these
surfaces.  For example, the latter authors proved that the only immersed
$H$-surfaces in $\mathcal{H}$ which are homeomorphic to the 2-sphere are
precisely the rotational $H$-spheres.  We mention that the classical
Alexandrov Theorem is not yet known in $\He$: ``Is a compact embedded
$H$-surface a rotational sphere".

It is natural (and we believe important) to solve the Dirichlet problem in $\He$; we do this here.

\section {Preliminaries}

\subsection{The Heisenberg space}

Let $\He$ denote the three-dimensional Heisenberg Lie group  endowed with a left invariant metric.  In fact, we have a  one-parameter family of metrics indexed by bundle curvature by a real parameter $\tau\neq 0$.  The spaces are simply connected  homogeneous Riemannian manifolds carrying a \mbox{$4$-dimensional} isometry group. In global
exponential coordinates they are $\R^3$ endowed in standard coordinates with the metrics
$$
ds^2=dx^2+dy^2+(\tau(ydx-xdy)+dz)^2.
$$
 A global orthonormal tangent frame is given by
$$
E_1=\partial_x-\tau y\partial_z,\;\;\;E_2=\partial_y+\tau x\partial_z,\;\;\;E_3=\partial_z.
$$
The corresponding Riemannian connection is
$\nab_{E_j}E_j=0$, $1\leqslant j\leqslant 3$, and
$$
\nab_{E_1}E_3=\nab_{E_3}E_1=-\tau E_2,\;\;\;\;
\nab_{E_2}E_3=\nab_{E_3}E_2=\tau E_1
$$
$$
\nab_{E_1}E_2=-\nab_{E_2}E_1=\tau E_3.
$$
In particular,
$$
[E_1,E_2]=2\tau E_3\;\;\;\mbox{and}\;\;\;  [E_1,E_3]=0=[E_2,E_3].
$$

The Heisenberg space is a Riemannian submersion $\pi\colon\,\He\to\R^2$
over the standard flat Euclidean plane $\R^2$
whose fibers are the vertical lines. Thus the fibers
are  the trajectories of a unit Killing vector field
and hence geodesics.
The horizontal vector fields $E_1,E_2$ are basic  since
they are the horizontal lifts of the vector fields of the orthonormal coordinate base of $\R^2$, namely, $\pi_*(E_1)=\partial_x$ and $\pi_*(E_2)=\partial_y$.

The isometries of the space are the translations generated by the Killing vector fields
$$
F_1=\partial_x+\tau y\partial_z,\;\;\;F_2=\partial_y-\tau x\partial_z,\;\;\;F_3=\partial_z,
$$
and the rotations about the $z$-axis corresponding to
$$
F_4=-y\partial_x+x\partial_y.
$$
The translations corresponding to $F_1$ and $F_2$ are,
respectively,
$$
(x,y,z)\mapsto (x+t,y,z+\tau ty)
$$
and
$$
(x,y,z)\mapsto (x,y+t,z-\tau tx)
$$
where $t\in\R$. Thus, by the group of isometries vertical planes go to vertical planes, and Euclidean lines
go to Euclidean lines. For additional information, we refer to \cite{da}.

\subsection{Graphs}

We denote by $S_0\subset\He$ the surface whose points satisfy $z=0$. Given a domain $\Omega\subset \R^2$ throughout the paper we also denote by $\Omega$ its lift to $S_0$.  We define the {\it graph\/} $\Sigma(u)$ of $u\in C^0(\bar\Omega)$ on $\Omega$ as
$$
\Sigma(u)=\{(x,y,u(x,y))\in{\cal H}:(x,y)\in\Omega\}.
$$
Consider the smooth function $u^*\colon\He\to\R$ defined as
$u^*(x,y,z)=u(x,y)$ and set $F(x,y,z)=z-u^*(x,y,z)$. Then $\Sigma(u)=F^{-1}(0)$, and therefore
$$
2H=\dive\left(\frac{\bar{\nabla} F}{|\bar{\nabla} F|}\right).
$$
Here $\dive$ and $\bar{\nabla}$ denote the divergence and gradient in $\He$ and the mean curvature function $H$ of the graph is with respect to the downward pointing normal vector.

We have
$$
\bar{\nabla} F=-(\tau y+u_x)E_1+(\tau x-u_y)E_2 + E_3.
$$
Since $E_1,E_2$ are basic, using the Riemannian submersion one shows that the $H$-graph  equation is
\be\label{pde}
\dive_{\R^2} \left(\frac{\alpha}{W}\partial_x
+ \frac{\beta}{W}\partial_y\right)+2H=0
\ee
where
$$
\alpha=\tau y+u_x,\;\;\;\;
\beta=-\tau x+u_y
$$
and
$$
W^2=1+ \alpha^2 + \beta^2.
$$
It follows easily that  $\Sigma(u)$ has mean curvature function $H$ if and only if $u$ is a solution of the following PDE
\be\label{edp}
Q_H(u):=\frac{1}{W^3}
\left((1+\beta^2)u_{xx} + (1+\alpha^2)u_{yy} -2\alpha\beta u_{xy}\right)+2H=0
\ee
for $\alpha$, $\beta$  and $W$ as above. We remark that this is the Euclidean mean curvature equation for $\tau=0$.

\subsection{Cylinders and cones}

Let $\gamma\colon\,I\to S_0\subset\He$ be a smooth curve parametrized on an interval $I\subset\R$ where $S_0$ is as above. We assume that $\gamma=\gamma(s)$ is parametrized so that $\bar \gamma=\pi\circ \gamma$ carries a parametrization by arc length. Thus  $\gamma(s)=(x(s),y(s),0)$ satisfies
$(x')^2+ (y')^2=1.$

The {\it  vertical cylinder\/} $\C_\gamma\subset\He$ over $\gamma$ is the surface generated by taking through each point of $\gamma$ the vertical geodesic fiber.  Thus $\C_\gamma$ is parametrized by $\varphi\colon\,I\times\R\to\He$
given by
$$
\varphi(s,t)=(x(s),y(s),t).
$$
Then, the mean curvature $H_\C$ (taken to be non-negative) of $\C_\gamma$ is
\be\label{h}
H_\C(s)=H_\C(s,t)=\frac{k(s)}{2}
\ee
where $k(s)$ is the geodesic curvature function of
$\bar \gamma$ with respect to the Euclidean metric.
Notice that $H_\C$ is independent of the parameter $\tau$.
To see that (\ref{h}) holds, first observe that the horizontal lift $T$ of $\bar{\gamma}\,'=d\bar \gamma/ds$ to each point of
$\C_\gamma$ forms a horizontal unit tangent vector field.  Since $\C_\gamma$ is ruled by vertical geodesics, it follows that the mean curvature of $\C_\gamma$ is
$2H_\C=\<\nab _TT, N\>$, where $N$ is the Gauss map of the cylinder $\C_\gamma$ chosen so that $H_\C$ is non-negative. But $N$ is the horizontal lift of a unit  normal vector field $\eta$ to $\bar\gamma$ in $\R^2$, and hence
$\<\nab _TT, N\>
=\<D_{\bar \gamma'}\bar \gamma',\eta\>=k$, where $D$ denotes the Euclidean connection.

 The {\it cone\/} ${\cal C}_\gamma\subset\He$ with vertex $P\in\He\backslash S_0$ and base curve $\gamma$ as above is just the Euclidean cone in $\R^3$ constituted of straight lines from $P$ through points of $\gamma$.
Thus ${\cal C}_\gamma$ is parametrized by
$$
\psi(s,t)=(1-t)P + t\gamma(s)
$$
where $t\in (0,+\infty)$.

Vertical lines remain invariant under the isometries of  $\He$. Thus the same holds for vertical  cylinders.
Also Euclidean lines are sent to Euclidean lines by isometries of $\He$, and vertical planes as well. Thus cones are also invariant by isometries.
Hence, to analyze the behavior of the mean curvature of a cone we may assume that the vertex is $P=(0,0,c)$ where $c\neq 0$. Then, either a computation using (\ref{edp}) or by a direct computation, we obtain that the mean curvature $H=H(s,t)$ of ${\cal C}_\gamma$ pointing down is given by
$$
H\!=\!\frac{ct^2(x^2+y^2+c^2)
(y''x'-x''y')}{2(\tau^2t^4(x^2+y^2)(x'y-y'x)^2
\!+\!2c\tau t^3(xx'+yy')(x'y-y'x)\!
+\!t^2(c^2+(x'y-y'x)^2 )^{3/2}}.
$$
Here the sign of $H$ is non-negative when $\gamma$ is a convex Jordan curve in $\R^2$.
In particular,
$$
H(s,1)\to H_\C(s)\;\;\;\mbox{as}\;\;\;c\to +\infty.
$$
and
$$
H(s_0,t)\to +\infty\;\;\;\mbox{as}\;\;\; t\to 0^+
$$
if $y''x'-x''y'>0$ at $\gamma(s_0)$.

We also have fixing  $t=t_0$ and allowing  $c\to+\infty$ that
$$
2H(s_0,t_0)\to (y''x'-x''y')(s_0),
$$
and this is also a proof that the mean curvature of a cylinder is given by (\ref{h}).

\section {The main result}

We now state and prove the Dirichlet theorem in Heisenberg space $\He$.

\begin{theorem}\po\label{ones} Let $\Omega\subset\R^2$ be a bounded domain with $C^3$ boundary $\Gamma=\partial\Omega$ whose curvature function  with respect to the inner orientation is $k>0$. Let $H$ be a constant satisfying  $0\leq 2H<k$ and let $\varphi\in C^0(\Gamma)$ be given.  Then there exists a smooth function $u$
satisfying  $u|_{\Gamma}=\varphi$ whose graph $\Sigma(u)$ in $\He$ has constant mean curvature $H$.

Moreover, if $M$ is a compact embedded connected surface inside the vertical cylinder $\C_\Gamma$ over $\Gamma$ with constant mean curvature  $H$, $\partial M=\partial \Sigma(u)$ and the mean curvature vector of $M$ points down, then  $M=\Sigma(u)$.
\end{theorem}

\proof  First suppose that $H=0$. In this case we prove a more general existence result. In fact, we allow
$k\geq 0$, and $\varphi$ to have a finite number of discontinuities
$E\subset\Gamma$, and at each discontinuity, $\varphi$
has a left and right limit.  The Nitsche graph (see \cite{ro}) $\gamma$
of $\varphi$ is the graph of $\varphi$ on
$\Gamma\setminus E$ together with  the vertical segments over each point of $E$, joining the left and right limits
of $\varphi$ at this point. The Nitsche graph $\gamma$ is a Jordan curve on the vertical cylinder
$\C_\Gamma$ and its vertical projection to $\Gamma$ is a monotone (constant on the vertical  segments) map.

Since $\C_\Gamma$ is mean convex with respect to the inside of $\C_\Gamma$, there is a least area embedded minimal disk $\Sigma$ inside $\C_\Gamma$ with $\partial\Sigma=\gamma$.

We claim that $\Sigma$ is a $z$-graph over $\Omega$ and solves the Dirichlet problem as desired. First observe that
$\Sigma$ is nowhere vertical.  To see this, suppose $p\in \mbox{int}\; \Sigma$ and the tangent plane to $\Sigma$ at $p$ is vertical. Let $\beta\in\R^2$
be a line such that the vertical plane $P=\pi^{-1}(\beta)$ equals the tangent plane to $\Sigma$ at $p$. Then $P\cap\Sigma$ near $p$ is an analytic curve topologically equivalent to $Re(z^k)$, $k\ge 2$, in a neighborhood of $z=0$. Each branch of these curves leaving $p$ must go to $P\cap\partial\Sigma=P\cap\gamma$, by the maximum principle, i.e., a cycle in $(\mbox{int}\;\Sigma)\cap P$ would bound a disk in $\Sigma$ and we could touch this disk at an interior point with
another vertical plane (which is  also a minimal surface).
Now $P\cap\gamma$ consists of two points of $\Gamma$, or
one or two vertical segments of $\gamma$, by convexity of
$\Gamma$. Hence, at least two of the branches of $P\cap\Sigma$ leaving $p$, go to the same point, or vertical segments of $\gamma$. This yields a compact cycle
$C\subset P\cap\Sigma$.  $\Sigma$ is simply connected so $C$ bounds a disk $D\subset\Sigma$. Using vertical planes in $\He$, we can touch $D$ at an interior point so $D$ would equal this vertical plane; a contradiction. Thus $\Sigma$
is nowhere vertical in its interior.

Now $\Sigma$ separates the vertical cylinder over $\Gamma$
into two components. So $\Sigma$ can be oriented with the unit normal pointing up in its interior. Then each vertical line over a point in the interior of $\Omega$, intersects $\Sigma$ in exactly one point, since at two successive points of intersection the normal to $\Sigma$ would point up and down. This proves $\Sigma$ is a graph over the interior of $\Omega$.

Now assume that $H\neq 0$ and $\varphi$ is continuous. We have seen that $u$ must be a solution of the
Dirichlet problem
\be\label{dir}
\hspace*{-65ex}\left\{ \begin{array}{l}
\! Q_H(u)=0
\vspace{1ex}\\
\! u|_{\Gamma}=\varphi
\end{array} \right.
\ee
where $Q_H$ was given in (\ref{edp}). To prove the existence part of the theorem, we use the continuity method. We show that the subset
$$
Z:=\{t\in\lbrack0,1]:\exists\, u_{t}\in C^3(\Omega)
\mbox{ such that }
Q_{tH}(u_{t})=0 \mbox{ and } u_{t}|_{\Gamma}=t\varphi\}
$$
is nonempty, open and closed in $[0,1]$. We have that $Z$ is not empty
since $0\in Z$; $S_0$ is a minimal surface in $\He$.  Standard arguments from the theory of quasilinear elliptic PDE's presented in  \cite{gt} give that $Z$ is open (a
consequence of the implicit function theorem). Moreover,  any solution of $Q_{H}(u)=0$ is smooth in $\Omega$.
Finally, that $Z$ is closed follows from the theory in  \cite{gt} once we show that a priori height and gradient estimates exist.

  We have from (\ref{edp}) that any Euclidean plane in $\R^3$ is a minimal surface in  $\He$.  In particular, each leaf of the foliation  of isometric surfaces $z=z_0=$ constant is minimal and diffeomorphic to the base $\R^2$ by the projection of the Riemannian submersion.
It follows using the maximal principle that any solution
$u$ of (\ref{dir}) satisfies
$$
u\geq\min_{\partial\Omega}\varphi.
$$

Fix a point $(x_0,y_0,0)\in \Omega$. Given
$z_0\in\R$, we consider the cone $C(z_0)$ with vertex $P=(x_0,y_0,z_0)$ constituted of straight lines from $P$ through points of the graph of $\varphi$ over $\Gamma$. Then, the piece $C_\varphi(z_0)$ of $C(z_0)$ from $P$ to the graph of $\varphi$ is contained inside the vertical cylinder over $\Gamma$. Notice that $C(z_0)$ is the cone $C_{\hat\Gamma}(z_0)$ over $\hat\Gamma=C_\Gamma(z_0)\cap S_0$.
Clearly, by choosing $z_0$ such that $|z_0|$ is large
enough, the geodesic curvature of $\hat\Gamma$ with
respect to the Euclidean metric is positive.
In fact, the curve converges to $\Gamma$ as $|z_0|\mapsto\infty$.
Therefore, by our previous discussion on the mean curvature of vertical cylinders and cones we have that choosing $z_0$ large enough, say $z_0=z_{1}$, and $z_0$ small enough, say $z_0=z_{2}$, that  $C(z_1)$ has mean curvature strictly
larger than $H$ everywhere and $C(z_2)$ has negative mean curvature (this cone is going down).  By the maximum principle, they are upper and lower barriers for the CMC $H$-graph equation on $\Omega$. Thus $C(z_1)$ and the above remark concerning planes below the graph of $\varphi$ provides an a priori height estimate for any solution of the Dirichlet problem (\ref{dir}) depending only on $\Omega$, $H$ and $\varphi$, that is,
$$
|u|_0\leqslant C_0(\Omega, H, \varphi).
$$
Moreover, the cones also provide the following bound along $\Gamma$ for the norm of the Euclidean gradient of $u$
$$
|Du|=\sqrt{u_x^2+u_y^2}\leqslant C_1(\Omega, H, \varphi).
$$

\vspace{1ex}
The next result  uses techniques developed in \cite{cns} to show  that global estimates of the gradient reduces to the boundary estimates already obtained.

\begin{lemma}\po\label{gradglob}
Let $u\in C^{3}(\Omega)\cap C^{1}(\bar \Omega)$ be a solution of (\ref{dir}). Assume that $u$ is bounded in $\Omega$ and that $|Du|$ is bounded in $\Gamma$.
Then $|Du|$ is
bounded in $\Omega$ by a constant that depends only on $|u|_0$ and $\sup_{\,\Gamma}|D u|$.
\end{lemma}

\proof To estimate $|D u|=\sqrt{u_x^2+u_y^2}$ in the interior of
$\Omega$ it suffices to obtain an estimate for
$\omega=\sqrt{\alpha^2+\beta^2}\ e^{Au}$ for some positive constant $A$ to
be chosen later. If $\omega$ achieves its maximum on $\Gamma$ then we have
the desired bound. Otherwise, $\omega$ must reach its maximum at an
interior point $p_0=(x_0,y_0)$ in $\Omega$.

We may choose coordinates of the ambient space such that
$$
\beta(p_0)=-\tau x_0+ u_y(p_0)=0.
$$
We denote
$$
v=\alpha(p_0)=\tau y_0+ u_x(p_0).
$$
The function $\phi=\ln\omega=\ln\sqrt{\alpha^2+\beta^2}+Au$ also takes a
maximum at $p_0\in\Omega$. That $\phi_x(p_0)=0$ yields
\be\label{one}
u_{xx}(p_0)=-Avu_x(p_0),
\ee
and $\phi_y(p_0)=0$ gives
\be\label{two}
u_{xy}(p_0)=-\tau(Avx_0+1).
\ee
Moreover, from $\phi_{xx}(p_0)\leqslant 0$ we obtain
\be\label{three}
vu_{xxx}(p_0)\leqslant
A^2v^3u_x(p_0)+A^2v^2u^2_x(p_0)-\tau^2(Avx_0+2)^2,
\end{equation}
and $\phi_{yy}(p_0)\leqslant 0$ yields
\be\label{four}
vu_{xyy}(p_0)\leqslant -Av^2u_{yy}(p_0)+\tau^2A^2x_0^2v^2-u_{yy}^2(p_0).
\end{equation}
On the other hand, from (\ref{edp}) and (\ref{one}) we have
\be\label{five}
u_{yy}(p_0)=-2H(1+v^2)^{1/2}+\frac{Av}{1+v^2}u_x(p_0).
\ee
Taking the derivative of (\ref{edp}) with respect to $x$ and using (\ref{one}) and (\ref{two}) yields
$$
u_{xxx}+ (1+v^2)u_{xyy}-2Av^2u_xu_{yy}
-2\tau^2v(A^2x_0^2v^2+3Ax_0v+2)-6AHv^2(1+v^2)^{1/2}u_x=0
$$
at the point $p_0$. Multiplying the last equation by $v$ and
using (\ref{five}) and inequalities (\ref{three}) and (\ref{four}) we obtain, after a long computation, that
$$
\frac{(v-\tau y_0)^2}{1+v^2}+\tau^2x_0^2\leqslant\frac{1}{A^2}(AG_1(v)+G_2(v))
$$
where
$$
G_1(v)=\frac{2H\tau y_0(1+v^2)^{1/2}}{v}+\frac{P(v)}{v^4},\;\;\;
G_2(v)=-4H^2+\frac{Q(v)}{v^4}
$$
and $\lim_{v\to\infty}P(v)/v^4=0=\lim_{v\to\infty}Q(v)/v^4$.
Therefore,
$$
\lim_{v\to\infty}G_1(v)=2H\tau y_0\;\;\;\;\mbox{and}\;\;\;\;
\lim_{v\to\infty}G_2(v)=-4H^2<0.
$$
It follows that we can choose $A>0$ such that
$$
\frac{(v-\tau y_0)^2}{1+v^2}+\tau^2 x_0^2\leqslant\frac{1}{2}.
$$
This gives an upper bound for $v^2$, and hence for $\omega=\sqrt{\alpha^2+\beta^2}\ e^{Au}$.
This concludes the proof of the Lemma.\qed
\vspace{1ex}

Hence $Z$ is closed, and this concludes the proof of the existence part of the Theorem for $0\leq 2H<k$. Now we prove that the graph $\Sigma=\Sigma(u)$ is unique. Suppose that
$M$ is an embedded \mbox{$H$-surface} inside the vertical cylinder
$\C_\Gamma$ over $\Gamma$ with $\partial M=\partial\Sigma$.
Then $M$ separates $\C_\Gamma$ into two components and we
assume the mean curvature vector of $M$ points into the lower component.  When the mean curvature  vector points toward the upper component, our argument will show that $M$ equals the
graph of the function $u$, equal to $\varphi$ on $\Gamma$,
with mean curvature $H$ and mean curvature vector pointing
toward the upper component.

 The mean curvature of the vertical cylinder over $\Gamma$
is strictly larger than $H$ and the mean curvature vector points inside the cylinder so the interior of $M$ is disjoint from the cylinder by the comparison principle.

Denote by $\Sigma(t)$ the surface $\Sigma$ translated $t$ by the flow of the Killing field $\partial z$. Since $\partial\Sigma$ is a $z$-graph, we have $\partial\Sigma(t)\cap\Sigma(0)=\emptyset$;  $\partial\Sigma(0)=\partial\Sigma$. Since $M$ is compact there is a $T>0$ such that
$\Sigma(T)\cap M=\emptyset$.

Now lower $\Sigma(T)$ to $\Sigma$ by the flow $\partial z$, letting $t$
go from $T$ to $0$. The mean curvature of each  $\Sigma(t)$ points down, so there can be no first  contact of $\Sigma(t)$ with $M$  for $t>0$, by the maximum principle. Thus $M$ is below $\Sigma$. Now choose $T<0$ so that $\Sigma(t)\cap M=\emptyset$. Move $\Sigma(T)$ up to $\Sigma$ by the flow $\partial z$, letting $t$ go from $T$ to $0$. There can be no first  contact of $\Sigma(t)$ with $M$ for $t\neq 0$ by the maximum principle (the mean curvature vector of $M$ points toward the downward component). Therefore $M$ is above $\Sigma$, and we obtain that $M=\Sigma$. This concludes the proof of the Theorem.\qed

\section {A further result}

It would be interesting to know if
Theorem \ref{ones} holds when we allow $2H=k$. In this section we give the following partial answer.

\begin{theorem}\po\label{second} Let $\Omega\subset\R^2$ be a bounded domain with $C^3$ boundary $\Gamma=\partial\Omega$ whose curvature function  with respect to the inner orientation is $k>0$. Let $H$ be a constant satisfying  $|\tau|/\sqrt{3}<H\leq k/2$. Then there exists a smooth function $u$ satisfying  $u|_{\Gamma}=0$ whose graph $\Sigma(u)$ in $\He$ has constant mean curvature~$H$.
\end{theorem}

We need  a supersolution $w$ defined in a neighborhood of $\Gamma$ (better than the cones in the preceding section);
$w$ is constructed in the following result.

\begin{proposition}\po\label{gradbound} Assume that
$u\in C^{2}(\Omega)\cap C^{1}(\bar \Omega)$ satisfies
$Q_{H}(u)=0$ in $\Omega$, $u|_{\Gamma}=0$ and
$|u|_0< M$. If $0<2H\le k$ on $\Gamma$, then there is a constant $C=C(H,\Omega, M)$ such that
$$
\sup_{\,\Gamma}|D u|\leq C.
$$
\end{proposition}

\proof Let $\gamma\colon\,[0,\ell]\rightarrow\Gamma$ be a parametrization by arc length and let $\nu$ stand for the unit normal vector to $\Gamma$ pointing to $\Omega$. We parametrize a neighborhood $U$ of $\Gamma$ in $\Omega$ by
\be\label{not}
P=P(s,t)=\gamma(s)+t\nu(s)
\ee
for $(s,t)\in [0,\ell]\times\lbrack0,\epsilon]$, where $0<\epsilon <1/k(s)$.
We compute (\ref{pde}) on $U$ making use of the orthonormal frame
$$
P_{t}=\nu,\;\;\;\;\frac{1}{\phi}P_{s}=\gamma'
$$
where $\phi(s,t)=1-tk(s)>0$. Notice that (\ref{pde}) can be written as
$$
Q_H(u)=\dive_{\R^2} \left(\frac{Z}{\sqrt{1+|Z|^2}}\right)+2H=0
$$
where $Z(p)=-\tau Jp + D u(p)$
and $J$ is the standard complex structure in $\R^2$.
Then,
\be\label{six}
W^3Q_H(u)=
W^3\dive_{\R^2} \left(\frac{1}{W}Z\right)+2HW^3=
-\frac{1}{2}\<D W^2,Z\> + W^2\dive_{\R^2}Z+2HW^3,
\ee
where $W^2=1+|Z|^2$.

We compute $W^3Q_H(w)=0$ for $w=w(t)$ to be chosen. Then
$D w=w_tP_t$ and
\be\label{seven}
W^2=1+|Z|^2=w_t^2-2\theta w_t+A
\ee
where $\theta=\tau\<JP,P_t\>=\tau\<\gamma,\gamma'\>$ and $A=1+\tau^2|\gamma+t\nu|^2$. Moreover,
$$
\dive_{\R^2}Z=\Delta w=w_{tt}-k_tw_t
$$
where
$$
k_t(s)=\<D_{P_s/\phi}P_s/\phi,P_ t\>,
$$
and hence, $k_0(s)=k(s)$. Thus,
\be\label{eight}
W^2\Delta w=w_t^2w_{tt}-2\theta w_tw_{tt}-k_tw_t^3
+2\theta k_tw_t^2+Aw_{tt}-Ak_tw_t.
\ee
Moreover,
$$
D W^2=(2w_tw_{tt}-2\theta w_{tt} + A_t)P_t
+(-2\theta_sw_t+A_s)\phi^{-2} P_s.
$$
Using $JP_t=-\phi^{-1}P_s=-\gamma'$ and $\phi^{-1}JP_s=P_t=\nu$, it is easy to see that
\be\label{nine}
\frac{1}{2}\<D W^2,Z\>
=w_t^2w_{tt}-2\theta w_tw_{tt}+\theta^2w_{tt}+Bw_t+C
\ee
where the functions $B$ and $C$ are bounded on $U$ and do not depend on $w$ or any of its derivatives.
It follows from (\ref{six}), (\ref{seven}), (\ref{eight}) and (\ref{nine}) that
$$
W^3Q_H(w)=2H(w_t^2-2\theta w_t+A)^{3/2}-k_tw_t^3
+2\theta k_tw_t^2+(A-\theta^2)w_{tt}-(Ak_t+B)w_t-C.
$$

For positive constants $L$ and $K$ choose
$$
w(t)=L\ln(1+K^{2}t).
$$
Then $w(0)=0$ and $w_{tt}=-w_{t}^{2}/L$.
Given $M>0$ choose $L=M/\ln(1+K)$.
Thus,
$$
w(t)=\frac{M}{\ln(1+K)}\ln(1+K^{2}t).
$$
Hence,
$$
w(1/K)=M
$$
and
$$
w_{t}(0)=\frac{MK^{2}}{\ln(1+K)}.
$$

We claim that we can choose $K>1/\epsilon$ large enough such that $Q_{H}(w)<0$ for all
$(s,t)\in [0,\ell]\times\lbrack0,1/K]$. This fact, together with $w(1/K)=M$ (recall that $|u|_0< M$) allows us to use $w$ as a barrier from above for $Q_{H}$ and conclude the proof.

It suffices to show
that $Q_{H}(w)<0$ at $t=0$ for $K$ large enough. Since
$w_{t}(0)\rightarrow +\infty$ as $K\rightarrow +\infty$, the
claim is clear at points of $\Gamma$ where $2H<k$.  If
$2H=k$ first observe that at $t=0$ 
$$
\lim_{K\rightarrow +\infty} \frac{(w_t^2-2\theta w_t+A)^{3/2}-w_t^3+2\theta w_t^2}{w_t^2}=-\theta.
$$
Then, we have that
$$
(A-\theta^2)w_{tt}(0)=-\frac{1}{L}(1+\tau^2(|\gamma|^2
-\<\gamma,\gamma'\>^2))w_t^2(0)<0,
$$
and the claim follows from the fact that $L\rightarrow 0^+$ as $K\rightarrow +\infty$.
\qed\vspace{2ex}

\noindent{\it Proof of Theorem \ref{second}:}
Let $\Omega(n)$ be the domain with boundary
$$
P(s,1/n)=\gamma(s)+\frac{1}{n}\nu(s)
$$
for large $n$, so $\partial\Omega(n)$ is smooth.
By Theorem \ref{ones} there exists an $H$-graph $\Sigma(n)$
with $\partial\Sigma(n)=\partial\Omega(n)$, since the curvature of $\partial\Sigma(n)$ is strictly greater than $2H$. Let $u_n$ be the function with graph $\Sigma(n)$

The curvature tensor of $\He$ is given for any $X,Y,Z\in T\He$ by
$$
R(X,Y)Z=-3\tau^2(X\wedge Y)Z
+ 4\tau^2R_1(\partial_z;X,Y)Z
$$
 where
$$
R_1(\partial_z;X,Y)Z=
\<Y,Z\>\<X,\partial_z\>\partial_z
+\<Y,\partial_z\>\<Z,\partial_z\>X
-\<X,Z\>\<Y,\partial_z\>\partial_z - \<X,\partial_z\>\<Z,\partial_z\>Y.
$$
Thus the (not normalized) scalar curvature of $\He$ is $S=-\tau^2$.

By Theorem 1 of \cite{ro2}, there is a positive constant $L$
such that $|u_n|_0\le L$ for each $n$. By the maximum principle, $u_{n+1}> u_n$ on the domain of $u_n$. Since the $u_n$ are uniformly bounded by $L$, the function
$$
u(x)=\lim_{n\to\infty}u_n(x),
$$
is well defined for $x\in\Omega$ and is an $H$-graph in $\Omega$. Moreover, the upper barrier $w$ constructed in Proposition \ref{gradbound} shows that $u$ takes the value zero on the boundary.\qed

{\renewcommand{\baselinestretch}{1}
\hspace*{-20ex}\begin{tabbing}
\indent \=Luis J. Alias \\
\>Departamento de Matematicas \\
\>Universidad de Murcia \\
\>  Campus de Espinardo E-30100 -- Spain\\
\> ljalias@um.es
\end{tabbing}}

{\renewcommand{\baselinestretch}{1}
\hspace*{-20ex}\begin{tabbing}
\indent \= Marcos Dajczer\\
\> IMPA \\
\> Estrada Dona Castorina, 110\\
\> 22460-320 -- Rio de Janeiro -- Brazil\\
\> marcos@impa.br\\
\end{tabbing}}

\vspace*{-2ex}

{\renewcommand{\baselinestretch}{1}
\hspace*{-20ex}\begin{tabbing}
\indent \= Harold Rosenberg\\
\> Departement de Mathematiques, \\
\>Universite de Paris VII,\\
\> 2 place Jussieu, 75251 -- Paris -- France\\
\> rosen@math.jussieu.fr
\end{tabbing}}

\end{document}